\documentclass[runningheads]{llncs}
\usepackage{graphicx}
\usepackage{mathtools}
\usepackage{amsmath}
\usepackage{amssymb}
\usepackage{scalerel,stackengine}

\newtheorem{observation}[remark]{Remark}

\def\one{\mbox{1\hspace{-4.25pt}\fontsize{12}{14.4}\selectfont\textrm{1}}}

\begin{document}
\title{Distribution of the sequence $[m]P$ in Elliptic Curves} 
%
%
\author{Markos Karameris\orcidID{0000-0001-7009-5285} 
}
\institute{National Technical University of Athens}
\maketitle              
\begin{abstract}
Major controversy surrounds the use of Elliptic Curves in finite fields as Random Number Generators. There is little information however concerning the "randomness" of different procedures on Elliptic Curves defined over fields of characteristic $0$. The aim of this paper is to investigate the behaviour of the sequence $\psi_m=[m]P$ and then generalize to polynomial seuences of the form $\phi_m=[p(m)]P$. We first study the sequence in the space of Elliptic Curves $E(\mathbb{C})$ defined over the complex numbers and then reconsider our approach to tackle real valued Elliptic Curves. In the process we obtain the measure with respect to which the sequence $\psi$ is equidistributed in $E(\mathbb{R})$. In Section \ref{sec} we prove that every sequence of points $P_n=(x_n,y_n,1)$ equidistributed w.r.t. that measure is not equidistributed$\mod(1)$ with the obvious map $x_n\to\{x_n\}$. 

\end{abstract}
\section*{Notation}

$\mathbb{Q}$: The field of rational numbers\\
$\mathbb{R}$: The field of real numbers\\
$\mathbb{C}$: The field of complex numbers\\
$\Lambda$: A complex lattice $\mathbb{Z}\omega_1\times\mathbb{Z}\omega_2$\\
$E(K)$: An Elliptic Curve defined over a subfield of the closed field $\bar{K}$\\
$\mathcal{C}(X)$: The algebra of continuous functions $X\to\mathbb{R}$\\
$\mathcal{R}(X)$: The algebra of Riemann integrable functions $X\to\mathbb{R}$\\
$g(\Lambda)$: The embedding of $\Lambda$ in the real plane\\
$\Delta$: The discriminant of an Elliptic Curve\\
$\mathcal{B}$: The Borel $\sigma-$algebra over a set $X$\\
$\mu$: A Borel measure over the corresponding algebra\\
$\wp$: The Weierstrass Elliptic Function on a lattice

\section{Introduction}
An elliptic curve is defined as a projective plane curve of genus $1$. It is a straightforward application of the Riemann-Roch theorem to obtain an equivalent Weierstrass equation of the curve $y^2=x^3+Ax+B$. The most important thing about Elliptic Curves that makes them interesting is the group structure we can endow them with. Thus performing the operation $mP$ for a point of the curve $P$ we get a new point on the curve. It is then natural to ask: How are these points distributed across the curve? Do we have an explosion towards infinity for example, with greater and greater leaps being made? To answer this question we will first examine the structure on an elliptic curve defined over $\mathbb{C}$.

\subsection{Elliptic Curves over $\mathbb{C}$}
An elliptic curve over $\mathbb{C}$ is actually isomorphic to a lattice over the complex numbers $\mathbb{C}/\Lambda$ where $\Lambda=\mathbb{Z}\omega_1\times\mathbb{Z}\omega_2$ with $\omega_{1,2}\in{\mathbb{C}}$. We also define the fundamental parallelogram as $D=\{a+x_1\omega_1+x_2\omega_2,  x_{1,2}<1 \text{ and } a\in\Lambda\}$. This isomorphism is provided by the Weierstrass function $\wp(z)$. The exact form of the isomorphism is in fact: $\phi:\mathbb{C}/\Lambda\to{}E(\mathbb{C}), z \longrightarrow{} (\wp(z),\wp'(z),1)$ and it is an isomorphism of Riemann surfaces.
In this context an isogeny between Elliptic Curves $E_1, E_2$ has the form of a map $\phi: \Lambda_1\to\Lambda_2$. The isogenies are actually exactly the maps of the form $\phi_a(z)=az\mod\Lambda_2$ where $a\in\mathbb{C}: a\Lambda_1\subseteq\Lambda_2$. In this context, an endomorphism of $E$ has the form $\phi(z)=az, a\Lambda\subseteq\Lambda$. Since each lattice corresponds uniquely to an elliptic curve, we can associate the $j-$invariant of the curve with the lattice as $j(\Lambda)$. Two Elliptic Curves are isomorphic iff $j(\Lambda_1)=j(\Lambda_2)$ or iff $a\Lambda_1=\Lambda_2$ for some $a\in\mathbb{C}$.
\begin{observation}
Suppose that $\omega_{1,2}$ is a basis for the lattice $\Lambda$. Then $n\omega_1+m\omega_2=\omega_2(\frac{\omega_1}{\omega_2}n+m)$ and thus $\Lambda=\omega_2\Lambda_{\tau}$ where $\Lambda_{\tau}=[\tau,1], \tau=\omega_1/\omega_2$. Thus every lattice can be written in the form $\Lambda_{\tau}, Im(\tau)>0$
\end{observation}

\section{Distribution in $E(\mathbb{C})$}
Since we will be studying functions that are periodic in a lattice it is essential to identify these functions and their behaviour.

\subsection{Fourier Series in Lattices}
\begin{observation}
Let $\Lambda$ be a real lattice $\Lambda=[1,\tau]$ and let $\tau_x,\tau_y$ be the projections of $\tau$ on the canonical vectors of $\mathbb{R}^2$. Then every function $f:\mathbb{R}^2/\Lambda\to{}A$ is double periodic in $\mathbb{R}^2$, or equivalently it can be identified with a function $f:\mathbb{R}^2\to{}A$ such that $\forall{}(x,y)\in\mathbb{R}^2, f(x,y)=f(x+1,y)=f(x+\tau_x,y+\tau_y)$.
\end{observation}
\begin{theorem}\label{th1}
Every function $f\in\mathcal{C}(\mathbb{R}^2/\Lambda)$ with $\Lambda=[1,\tau]$ admits a Fourier series expansion of the form:
\begin{equation}
f(x,y)=\sum_{n_1,n_2\in\mathbb{Z}}a_{n_1,n_2}e^{2\pi{}j(n_1x+\frac{n_2-n_1\tau_x}{\tau_y}y)}
\end{equation}
\end{theorem}

\begin{lemma} \label{lem1}
Define the transformation $A= \left[\begin{matrix} 1 & 0 \\ \tau_x & \tau_y \end{matrix}\right]$. Then $A$ maps $\mathbb{R}^2/[0,1]\times{}[0,1]$ to $\mathbb{R}^2/\Lambda$ continuously. (By the same methods we can also prove the continuity of $A^{-1}=\left[\begin{matrix} 1 & 0 \\ -\frac{\tau_x}{\tau_y} & \frac{1}{\tau_y} \end{matrix}\right]$ since they have the exact same form)
\begin{proof}
For every pair of points: $|A[x_1,y_1]-A[x_2,y_2]|=|[(x_1-x_2)+\tau_x(y_1-y_2),\tau_y(y_1-y_2)]|$, setting $x_1-x_2=x', y_1-y_2=y'$ we obtain: $|A[x_1,y_1]-A[x_2,y_2]|=|[x'+\tau_xy',\tau_yy']|={x'}^2+2\tau_xy'x'+{\tau_y}^2{y'}^2+{\tau_y}^2{y'}^2$ using the Cauchy-Schwartz inequality: $|A[x_1,y_1]-A[x_2,y_2]|\leq{}\max({x'}^2+{y'}^2+{y'}^2+{\tau_x}^2{x'}^2,|{\tau}|^2({x'}^2+{y'}^2)+{y'}^2+{\tau_x}^2{x'}^2)$ if $|{\tau}|^2>1$ or $|{\tau}|^2\leq{}1$. With the exact same logic for ${\tau_x}^2$ we get that $|A[x_1,y_1]-A[x_2,y_2]|\leq{}|[x_1-x_2,y_1-y_2]|(2+|{\tau}|^2)$. We have thus shown uniform continuity.
\end{proof}
\end{lemma}

\begin{theorem}
$f\in\mathcal{C}(\mathbb{R}^2/[0,1]\times[0,1])$ iff $f\circ{}A^{-1}\in\mathcal{C}(\mathbb{R}^2/\Lambda)$.
\begin{proof}
Suppose $f(x,y)=f(x+1,y)=f(x,y+1)$ then $f\circ{}A^{-1}(x+1,y)=f(A^{-1}[x+1,y])=f(x+1+\frac{\tau_x}{\tau_y}y,\frac{1}{\tau_y}y)=f(x+\frac{\tau_x}{\tau_y}y,\frac{1}{\tau_y}y)=f\circ{}A^{-1}(x,y)$ and $f\circ{}A^{-1}(x+\tau_x,y+\tau_y)=f(x+\tau_x-\frac{\tau_x}{\tau_y}y-\tau_x,\frac{1}{\tau_y}y+1)=f(x+\frac{\tau_x}{\tau_y}y,\frac{1}{\tau_y}y)=f\circ{}A^{-1}(x,y)$.
For the inverse it suffices to assume $f\in\mathcal{C}(\mathbb{R}^2/\Lambda)$ and then we have $f\circ{}A(x+1,y)=f(A[x+1,y])=(x+\tau_xy+1,\tau_yy)=f(A[x,y])=f\circ{}A(x,y)$ and $f\circ{}A(x,y+1)=f(A[x,y+1])=f(x+\tau_xy+\tau_x,\tau_yy+\tau_y)=f(x+\tau_xy,\tau_yy)=f(A[x,y])=f\circ{}A(x,y)$. The continuity of each of these composite functions follows from Lemma \ref{lem1}.
\end{proof}
\end{theorem}
\begin{corollary}
For each $f\in\mathcal{C}(\mathbb{R}^2/\Lambda)$ there is exactly one coresponding $f\in\mathcal{C}(\mathbb{R}^2/[0,1]\times[0,1])$.
\end{corollary}
We now finish the proof of theorem \ref{th1}:
\begin{proof}
Suppose $f\in\mathcal{C}(\mathbb{R}^2/\Lambda)$ then we define $A$ as before according to the values of the lattice $\Lambda$. We now get a function $f\circ{}A\in{}\mathcal{C}(\mathbb{R}^2/[0,1]\times[0,1])$ and thus $f\circ{}A$ admits a Fourier series expression of the form $f(A[x,y])=\sum_{n_1,n_2\in\mathbb{Z}}e^{2\pi{}j[n_1,n_2][x,y]}$. The Fourier series expression of $f$ is then $f(x,y)=\sum_{n_1,n_2\in\mathbb{Z}}e^{2\pi{}j[n_1,n_2]A[x,y]}=\sum_{n_1,n_2\in\mathbb{Z}}e^{2\pi{}j[n_1,n_2][x-\frac{\tau_x}{\tau_y}y,\frac{1}{\tau_y}y]}=$\\
$\sum_{n_1,n_2\in\mathbb{Z}}a_{n_1,n_2}e^{2\pi{}j(n_1x+\frac{n_2-n_1\tau_x}{\tau_y}y)}$.
\end{proof}
\begin{observation}
In this section we only worked with lattices of the form $[\tau,1]$ but it is possible to work with any two vectors $[u,v]$ defining a lattice (which means linearly independend). Then the general form of the Fourier transform is $f(x,y)=\sum_{n_1,n_2\in\mathbb{Z}}e^{2\pi{}j[n_1,n_2]A^{-1}[x,y]}$ where $A=\left[\begin{matrix} u_x & u_y \\ v_x & v_y \end{matrix}\right]$.
\end{observation}

This section aims to show one thing basically which is now immediate:
\begin{theorem} \label{thth2}
The sub-algebra of trigonometric polynomials with variables of the form $e^{2\pi{}j(n_1x+\frac{n_2-n_1\tau_x}{\tau_y}y)}$ is dense in $\mathcal{L}(\mathbb{R}^2/\Lambda)$.
\begin{proof}
Since $\mathcal{C}([a,b])$ is dense in $\mathcal{L}^2$ (w.r.t. the $sup$ metric) and trigonometric polynomials are dense in $\mathcal{C}([a,b])$ as a consequence of Theorem \ref{th1}, the result is immediate.
For a proof of the density of $\mathcal{C}([a,b])$ in $\mathcal{L}^p, p\geq{}1$ see \cite{realan} page 153.
\end{proof}
\end{theorem}

\subsection{Equidistribution of $[m]P$ in $\Lambda$}

Throughout this section we will be working with the map $g:\mathbb{C}\to\mathbb{R}^2$ sending $z_x+z_yi\to{}(z_x,z_y)$. This map sends $\Lambda$ to a real valued lattice in $\mathbb{R}^2$ and we can then define equidistribution in the usual way for a compact metric space.
\begin{definition}
A sequence $s_n$ in a compact metric space $X$ equiped with the Borel probibility measure $\mu$ is equidistributed if $\lim_{n\to\infty}\frac{1}{n}\sum_{i=0}^{n-1}f(s_i)=\int_{X}fd\mu$ for every Riemann integrable $f:X\to\mathbb{C}$.
\end{definition}
\begin{observation}
A sequence $z_n$ is equidistributed in $\mathbb{C}/\Lambda$ iff for every $f\in{}\mathcal{R}(\mathbb{R}^2/g(\Lambda)), f:\mathbb{R}^2/g(\Lambda)\to\mathbb{R}$ we have \\ $\lim_{n\to\infty}\frac{1}{n}\sum_{i=0}^{n-1}f(g(z_n))=\frac{1}{\mu_{\mathbb{R}}(g(\Lambda))}\int_{g(\Lambda)}f(x,y)dxdy$. The use of $dxdy$ instead of $d\mu$ follows from the function being Riemann Integrable.
\end{observation}
\begin{theorem}
A sequence $z_n$ is equidistributed in $\mathbb{C}/\Lambda$ iff
\begin{equation}
\lim_{N\to\infty}\frac{1}{N}\sum_{n=0}^{N-1}e^{2\pi{}j(n_1{z_n}_x+\frac{n_2-n_1\tau_x}{\tau_y}{z_n}_y)}=0, \forall{}n_1,n_2\in\mathbb{Z}
\end{equation}
\begin{proof}
($\implies$) This part is immediate since we just have to substitute $f(x,y)=e^{2\pi{}j(n_1x+\frac{n_2-n_1\tau_x}{\tau_y}y)}$.\\
($\impliedby$) From Theorem \ref{thth2} we can see that trigonometric polynomials are dense in $\mathcal{R}(\mathbb{R}^2/g(\Lambda))$. A standard limit argument similar to the $\mathbb{R}$ case now implies the result.
\end{proof}
\end{theorem}
\begin{theorem}
For a point $z$, the sequence $nz\mod\Lambda=z_{nx}+z_{ny}i$ is equidistributed in $\mathbb{C}/\Lambda$ iff $n_1{z_n}_x+\frac{n_2-n_1\tau_x}{\tau_y}{z_n}_y\not\in\mathbb{Z}$ for every choice of $(n_1,n_2)\ne{}(0,0)$.
\begin{proof}
Setting $k(n_1,n_2)=n_1z_x+\frac{n_2-n_1\tau_x}{\tau_y}z_y$ we get $\lim_{N}\frac{1}{N}\sum_{n=0}^{N-1}e^{2\pi{}jnu}$ and thus if $k(n_1,n_2)\in\mathbb{Z}$ for some $(n_1,n_2)\in\mathbb{Z}^2/(0,0)$ we have \\ $\lim_{N\to\infty}\frac{1}{N}\sum_{n=0}^{N-1}1=1\ne{}0$. Otherwise we have $|\frac{1}{N}\sum_{n=0}^{N-1}e^{2\pi{}jnu}|\leq\frac{1}{N}\frac{|e^{2\pi{}ju(N-1)}-1|}{|e^{2\pi{}ju}-1|}\leq\frac{1}{N}\frac{2}{|e^{2\pi{}ju}-1|}$ and thus $\lim_{N\to\infty}\frac{1}{N}\sum_{n=0}^{N-1}e^{2\pi{}jnu}|=0$.
\end{proof}
\end{theorem}

A few obvious families of points where equidistribution fails are points parallel to one of the lattice defining vectors:
\begin{enumerate}
    \item $\forall{}(x,y): y=0 \implies n_1x=0$ and thus a solution for $k(n_1,n_2)=0$ will always be $(0,n), \forall{}n\in\mathbb{Z}$.
    \item $\forall{}(x,y): (x,y)=(a\tau_x,a\tau_y)$ we have $k(n_1,n_2)= an_1\tau_x+n_2a-n_1a\tau_x=0 \implies n_2=0$  and thus we obtain a solution for $k(n_1,n_2)$ which is $(n,0),\forall{}n\in\mathbb{Z}$.
    \item all elements parallel to the diagonals: $\forall{}(x,y)=(\lambda\tau_x+\lambda,\lambda\tau_y)$ we have $k(n_1,n_2)=\lambda{}n_1\tau_x+\lambda(n_1+n_2)-\lambda{}n_1\tau_x=(n_1+n_2)\lambda$ and thus an obvious solution is $(n_1,n_2)=(n,-n), n\in\mathbb{Z}$.
\end{enumerate}

\section{Real Elliptic Curves}
So far we have studied the equidistribution in complex Elliptic Curves. We will now shift our focus to Elliptic Curves $E(\mathbb{R})$. Naturally we first study the values $z\in\mathbb{C}/\Lambda$ for which $\wp(z)\in\mathbb{R}$. A more detailed analysis with applications can be found in \cite{realpart}.
\subsection{The Real Part of $\wp$}
\begin{theorem}
Let $\Lambda$ correspond to the Elliptic Curve $y^2=4x^3+g_2x+g_3$ where $g_2=g_2(\Lambda),g_3=g_3(\Lambda)$ are the invariants of the lattice. Then $g_2,g_3\in\mathbb{R} \iff \Lambda$ is invariant under complex conjugation.
\begin{proof}
$(\impliedby)$ is obvious since $g_2(\Lambda)=\sum_{\omega\in\Lambda^*}\frac{1}{\omega^4}$ and $g_3(\Lambda)=\sum_{\omega\in\Lambda^*}\frac{1}{\omega^6}$ and thus $g_2=\bar{g_2}$ and $g_3=\bar{g_3}$.\\
($\implies$) We know $\wp(z)=\frac{1}{z^2}+\sum_{n=1}^{\infty}(2n+1)G_{2n+2}(\Lambda)z^{2n}$ where $G_{2n+2}(\Lambda)$ are the Eisenstein series of weight $2n+2$ of the lattice. Setting $a_1=g_2/20, a_2=g_3/28$ and $(2n+1)G_{2n+2}(\Lambda)=a_n$ in general we get: $\wp(z)=\frac{1}{z^2}+\sum_{n=1}^{\infty}a_nz^{2n}$. By differentiating the Weierstrass equation we get $\wp"(z)=6{\wp(z)}^2-\frac{g_2}{2}$. By comparing the coefficients of $z^{2n}$ we have:\\ $a_{n+1}=\frac{6}{(2n+1)(2n+2)-12}\sum_{i=1}^ka_ka_{n-k}$. Thus inductively we get that $a_n\in\mathbb{R},\forall{}n\in\mathbb{N}$ and thus $\bar{\wp(z)}=\wp(\bar{z})$. This implies $\frac{1}{\bar{z}^2}+\sum_{\omega\in\Lambda}\frac{1}{(\bar{z}+\omega)^2}-\frac{1}{\omega^2}=\frac{1}{\bar{z}^2}+\sum_{\omega\in\Lambda}\frac{1}{(\bar{z}+\bar{\omega})^2}-\frac{1}{{\bar{\omega}}^2}$ form which we finally have: $\omega\in\Lambda \iff \bar{\omega}\in\Lambda$.
\end{proof}
\end{theorem}

\begin{corollary}
If $x\in\mathbb{Z}$ then the above theorem implies that for any Elliptic Curve with $g_2,g_3\in\mathbb{R}$ we have $\bar{\wp(x)}=\wp(\bar{x})=\wp(x)$ and $\bar{\wp(jx)}=\wp(-jx)=\wp(jx)$ and thus all purely real and imaginary values are in $\mathbb{R}$.
\end{corollary}

Let ${\wp'(z)}^2=4(\wp(z)-e_1)(\wp(z)-e_2)(\wp(z)-e_3)$ and observe that $\wp(z)=e_i \iff \wp'(z)=0$ and this only happens in the half-periods of the lattice.
Now consider two cases:
\begin{itemize}
    \item if $\Delta={g_2}^3-27{g_3}^2>0$ then $e_i\in\mathbb{R}$ and setting $e_1>e_2>e_3$ we can write $\wp(\frac{\omega_i}{2})=e_i$ where $\omega_2=\omega_1+\omega_3$ and $\Lambda=[\omega_1,\omega_3]$. Taking into account the fact that $\wp$ assumes every value in $\mathbb{R}$ exactly twice in $[0,\omega_1],[0,\omega_3],[\omega_1,\omega_1+2\omega_3],[\omega_3,\omega_3+2\omega_1]$ we have the full set of points where $\wp$ is real. Note that we have a square lattice.
    \item if $\Delta={g_2}^3-27{g_3}^2<0$ then we have two complex roots $e_1,e_3$ and one real root $e_2$. Then we have $\wp(\frac{\omega_2}{2})=e_2$ with every real value attained exactly twice both on the real and imaginary axis \\ $[-\frac{\omega_2}{2},\frac{\omega_2}{2}],[\frac{\omega_1-\omega_3}{2},\frac{\omega_3-\omega_1}{2}]$. Note the rhombic shape of the lattice.
\end{itemize}

\begin{figure}[htp]
\centering
\includegraphics[width=10cm]{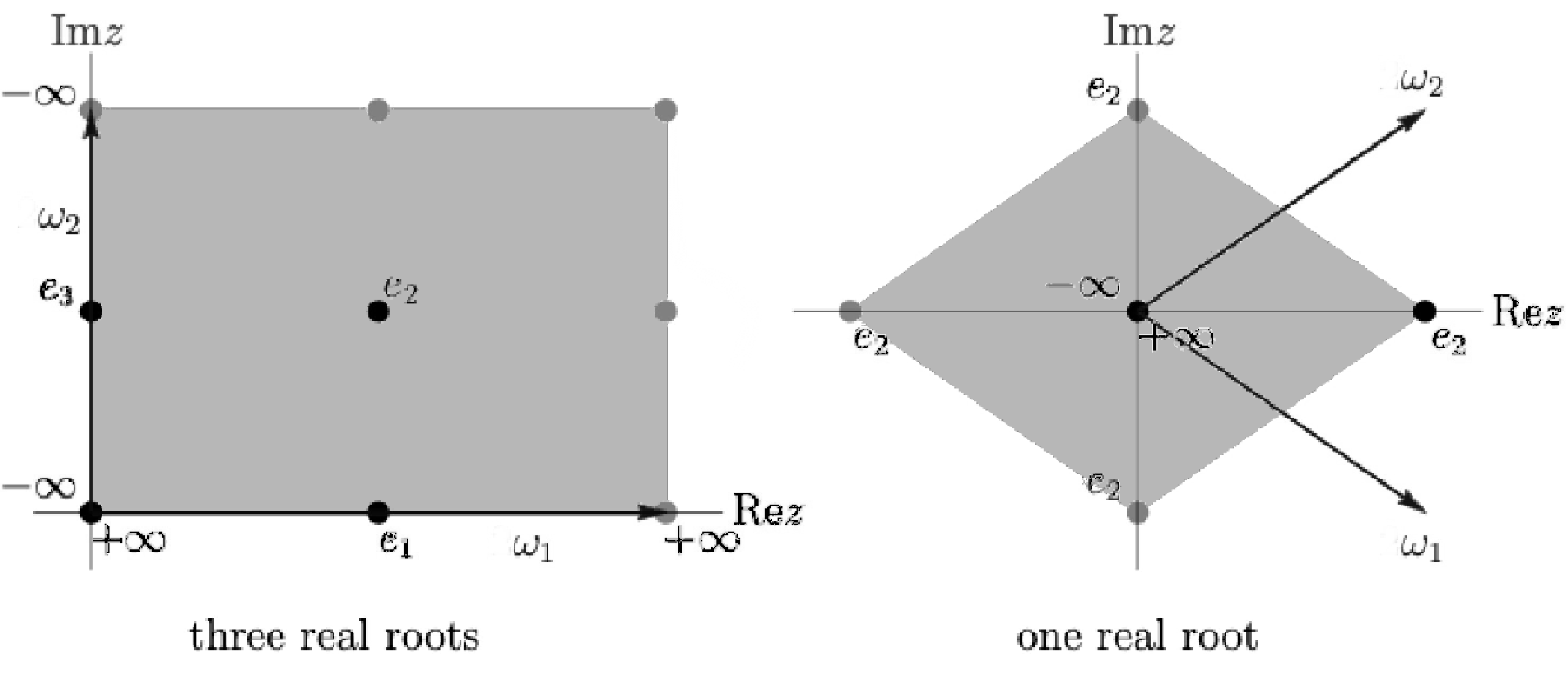}
\caption{The real part of $\wp$ when $\Delta>0$ and $\Delta<0$}
\label{fig:rom}
\end{figure}

\begin{observation}
Since we are only looking at real Elliptic Curves we only need to consider the values $x\in{}[e_3,e_2]\cup{}[e_1,+\infty]$ for the case with three real roots case. That is the intervals: $[0,\omega_1]$ and $[\frac{\omega_3}{2},\frac{\omega_3}{2}+\omega_1]$.\\
For the case of one real root we only need to consider the interval $[\frac{\omega_2}{2},\frac{\omega_2}{2}]$. Since $\wp$ is double periodic we can equivalently consider the set $[0,\omega_2]$ so that we have the same form in both cases.
\end{observation}

\begin{observation}
Note that every single point where $\wp$ is real valued is either parallel to the lattice vectors or on the diagonal. This means equidistribution fails for those points and indeed it should! The set $X=\{z\in\mathbb{C}/\Lambda: \wp(z)\in\mathbb{R}\}$ has measure $0$ in the probability space we defined previously. Also every sequence in $X$ will stay in the set (which has measure $0$) and thus there is no way it will exhibit the recurrence properties expected from equidistributed sequences.
\end{observation}

\subsection{Equidistribution in $E(\mathbb{R})$}
Let us begin by noting that since it is more convenient to deal with points on the real axis for $\wp$ we will keep the standard coordinates defined in the above section $\Lambda=[\omega_1,\omega_3]$. We will thus \textbf{not} transform the rhomboid lattice as usual by multiplying with $1/\omega_1$.
As noted in the previous section we also consider two cases here:
\begin{itemize}
    \item When $\Delta>0$ we look at the set where $y^2\geq{}0$, that is: $A_1=[0,\omega_1]\cup{}[\frac{\omega_3}{2},\omega_1+\frac{\omega_3}{2}]$ in which every value of $\wp(z)$ appears twice as $(\wp(z),\wp'(z),1)$ and $(\wp(z),-\wp'(z),1)$ in symmetric values of $z$ as $\wp(z)=\wp(-z)=\wp(\omega_1-z)$ and $\wp'(z)=-\wp(-z)=-\wp(\omega_1-z)$.
    \item When $\Delta<0$ we examine the set $A_2=[0,\omega_2]$. The same here is true for the values of $\wp$.
\end{itemize}
\begin{theorem} \label{ththe}
Define the probability space $(A_1,\mathcal{B}_1,\mu_+)$ where \\ $\mu_+(X)=\frac{\mu(X\cap{}[0,\omega_1])+\mu(X\cap[\frac{\omega_3}{2},\omega_1+\frac{\omega_3}{2}])}{2|\omega_1|}$ then the sequence $s_n=nz$ is equidistributed when $z\in\mathbb{R}$.
\begin{proof}
Let $X=[0,\omega_1]$, then $nz\in\mathbb{R}$ and thus $\lim_{N\to\infty}\frac{|\{s_1,\dots,s_n\}\cap{}A_1|}{N}=\mu(X)=\frac{\mu(A_1)}{2}$.
\end{proof}
\end{theorem}
The damage can be minimized by considering both of these probability spaces separately like so:
$B_1=[0,\omega_1]$ and $\Gamma_1=[\frac{\omega_3}{2},\omega_1+\frac{\omega_3}{2}]$ where both $B_1$ and $\Gamma_1$ are measure preserving systems under the transform $T(z)=z+a, a\in\mathbb{R}$. The first thing we observe is that in this case we have a space isomorphic to $[0,1]$ and thus we can use Weyl's Criterion.
\begin{definition}
 We say that a sequence $s_n\in{}[a,b]$ is equidistributed in $[a,b]$ iff $\lim_{N\to\infty}\frac{|\{s_1,\dots,s_n\}\cap[c,d]|}{N}=\frac{d-c}{b-a},$ for every $[c,d]\subseteq{}[a,b]$.
\end{definition}
\begin{theorem} (\textbf{Weyl's Criterion})\\
Suppose we have a sequence $s_n\in\mathbb{R}$, then the following are equivalent:
\begin{enumerate}
    \item $s$ is equidistributed in $[a,b]$
    \item for every $f$ Riemann integrable in $[0,1]$ it holds that\\ $\lim_{N\to\infty}\frac{1}{N}\sum_{n=1}^{N-1}f(s_n)=\frac{1}{b-a}\int_a^bf(x)dx$
    \item $\forall{}k\in\mathbb{Z}^*$ we have $\lim_{N\to\infty}\frac{1}{N}\sum_{n=1}^Ne^{\frac{2\pi{}jks_n}{b-a}}=0$
\end{enumerate}
\end{theorem}
For more details and applications on Weyls Criterion see \cite{ref_url1}.
\begin{observation}
If $s_n$ is equidistributed w.r.t. $\mu$ then $s_n+a$ is also equidistributed with respect to $\mu$ for every $a\in{}[0,\omega]$. This is immediate since $\lim_{N\to\infty}\frac{1}{N}\sum_{n=1}^Ne^{\frac{2\pi{}jk(s_n+a)}{b-a}}=e^{\frac{2\pi{}jka}{b-a}}\lim_{N\to\infty}\frac{1}{N}\sum_{n=1}^Ne^{\frac{2\pi{}jks_n}{b-a}}$ and thus it does not affect convergence to $0$. Indeed viewing $[0,\omega]$ as a topological group $\mathbb{R}/\omega\mathbb{Z}$ with addition, we get that $\mu$ is the normalized Haar measure as it is shift invariant, regular and suported on the whole $[0,\omega]$.
\end{observation}
\begin{lemma} \label{thislem}
For a point $z\in\mathbb{C}/\Lambda$ the sequence $s_n=nz$ is equidistributed in $(B_1,\mathcal{B}_1,\mu)$ (and $(A_2,\mathcal{B}_2,\mu)$ equivalently) iff $\lim_{N\to\infty}\frac{1}{N}\sum_{n=1}^{N-1}e^{\frac{2\pi{}jks_n}{\omega_1}}=0$ \\ (or $\lim_{N\to\infty}\frac{1}{N}\sum_{n=1}^{N-1}e^{\frac{2\pi{}jks_n}{\omega_2}}=0$ equivalently).
\begin{proof}
Given a space $[0,t]$ the Fourier expansion of fanctions $f\in{}C([0,t])$ is given by $f=\sum_{n=1}^{\infty}e^{2\pi{}jnx/t}$. The density of these trigonometric polynomials now follows and from the exact same argument in the proof of Weyl's Criterion we obtain $s_n$ is equidistributed in $[0,t]$ if and only if $\lim_{N\to\infty}\frac{1}{N}\sum_{n=1}^{N-1}e^{\frac{2\pi{}jks_n}{t}}=0$. Since in both cases $z,\omega_1,\omega_2\in\mathbb{R}$ the result immediately follows.
\end{proof}
\end{lemma}

\begin{theorem} \label{th2}
For a point $z\in\mathbb{C}/\Lambda$ the sequence $s_n=nz$ is equidistributed in $(B_1,\mathcal{B}_1,\mu)$ (and $(A_2,\mathcal{B}_2,\mu)$ equivalently) iff $z\not\in\omega_1\mathbb{Q}$ (or equivalently $z\not\in\omega_2\mathbb{Q}$). Then if $P=(\wp(z),\wp'(z))$ and $[n]P$ is equidistributed, we obtain that $P$ is not an element of the torsion subgroup of the curve $E_{Tor}=\{P\in{}E:[m]P=0,\text{for some }m\in\mathbb{Z}\}$.
\begin{proof}
By Lemma \ref{thislem} we have that $s_n$ is equidistributed in $(B_1,\mathcal{B}_1,\mu)$ iff \\ $\lim_{N\to\infty}\frac{1}{N}\sum_{n=1}^{N-1}e^{\frac{2\pi{}jknz}{\omega_1}}=0 \iff \lim_{N\to\infty}\frac{1}{N}\frac{(e^{\frac{2\pi{}jkz}{\omega_1}})^N-1}{e^{\frac{2\pi{}jkz}{\omega_1}}-1}=0, \forall{}k\in\mathbb{Z}^*$ where the last expression can only occur when $z\not\in\omega_1\mathbb{Q}$. Indeed if $z\not\in\omega_1\mathbb{Q}$ then we can choose $k\mathbb{Z}^*: \frac{kz}{\omega_1}\in\mathbb{Z}$ and thus $\lim_{N\to\infty}\frac{1}{N}\sum_{n=1}^{N-1}e^{\frac{2\pi{}jks_n}{\omega_1}}=1$. If $z\not\in\omega_1\mathbb{Q}$ then $\frac{(e^{\frac{2\pi{}jkz}{\omega_1}})^N-1}{e^{\frac{2\pi{}jkz}{\omega_1}}-1}\leq{}\frac{2}{e^{\frac{2\pi{}jkz}{\omega_1}}-1}$ and the result follows. The proof for $(A_2,\mathcal{B}_2,\mu)$ is the same.
\end{proof}
\end{theorem}

\begin{observation}
For every interval $(x_a,x_b)$ with $y>0$ or $y<0$ we have a unique interval $\wp^{-1}((x_a,x_b))\in{}[0,\omega]$ with $\mu(\wp^{-1}((x_a,x_b)))=|\wp(x_a)-\wp(x_b)|$. This shows that if $s_n$ is equidistributed in $[0,\omega]$ it is also dense in $[0,\omega]$ which means that $\wp(s_n)$ is dense in $[e,\infty)$ and thus $(\wp(s_n),\wp'(s_n))$ is dense in $E(\mathbb{R})$. In the case of $[n]P$ this means that either $[m]P=O$ for some $m\in\mathbb{Z}$ or $[m]P$ is dense in $E(\mathbb{R})$.
\end{observation}

Let us refer to both $\omega_1,\omega_2$ as $\omega$ for simplicity, since both cases yield the same result.
However in the case of $\omega_1$ the real period of the associated Elliptic Curve is actually $2\omega_1$ since we have two connected components but we ommited the $\Gamma_1$ so we proceed similarly. Basically we consider $\omega=\int_{e}^{\infty}\frac{dx}{y}$ where $e$ is $e_1$ or $e_2$ in each case.
Then returning to Weyl's Criterion we obtain the following result:

\begin{corollary} \label{thiscor}
Let $z\not\in\omega\mathbb{Q}$ then the sequence $z_n=nz$ is equidistributed in $[0,\omega]$ and thus $\forall{}f\in{}\mathcal{R}([0,1])$ we have $\lim_{N\to\infty}\frac{1}{N}\sum_{n=1}^{N}f(z_n)=\frac{1}{\omega}\int_0^{\omega}f(z)dz$.
\end{corollary}

\begin{observation}
In the above corollary we are only considering Riemann Integrable functions and so the use of the differential $dz$ is equivalent to using Lebesgue integreation w.r.t. $\mu$. Notice that the $\frac{1}{\omega}$ term appears since we are using the normalized measure $\mu_{\omega}=\mu/\omega$.
\end{observation}

Before moving to the main theorem we clarify the following:
\begin{definition}
 We say that a function $f:(a,\infty)\to\mathbb{R}$ is improper Riemann integrable and write $f\in\mathcal{IR}((a,\infty))$ iff $\lim_{{\substack{\epsilon\to{}0\\r\to\infty}}}\int_{a+\epsilon}^{r}f(x)dx=c\in\mathbb{R}$.
\end{definition}

Corollary \ref{thiscor} enables us to shift to points on the real curve:
\begin{theorem}
Let $z\not\in\omega\mathbb{Q}$, then the sequence $z_n=nz$ is equidistributed in $[0,\omega]$ and for every $f$ bounded in $[e,\infty)$ such that $\frac{f(x,\pm{}y)}{y}\in{}\mathcal{IR}((e,\infty))$,\\ $\lim_{N\to\infty}\frac{1}{N}\sum_{n=1}^{N}f(\wp(z_n),\wp'(z_n))=\frac{1}{\omega}\int_{e}^{\infty}(f(x,y)+f(x,-y))\frac{dx}{y}$ where $e=e_1$ or $e=e_2$ depending on the case of $\omega_{1,2}$ and $y\geq{}0$.\footnote[1]{ In this theorem $y$ is treated as a function of $x$ by seperating the parts $y>0$ and $y<0$ and thus $f$ is not a two variable function but rather a function of $x$ only.}
\begin{proof} One obvious obstacle is that $f(\wp(\omega),\wp'(\omega))$ is not defined since $\wp(\omega)$ is not defined in $\mathbb{R}$. We can fix that however by setting $f(\wp(\omega),\wp'(\omega))$ equal to any value or even better $f(\wp(\omega),\wp'(\omega))=\lim_{z\to\omega}f(\wp(z),\wp'(z))$ if it exists. Then from Corollary \ref{thiscor} it immediately follows that: \\
$\lim_{N\to\infty}\sum_{n=1}^{N}f(\wp(z_n),\wp'(z_n))=\int_0^{\omega}f(\wp(z),\wp'(z))\frac{dz}{\omega}= \\ \int_0^{\omega/2}f(\wp(z),\wp'(z))\frac{dz}{\omega}+\int_{\omega_2/2}^{\omega}f(\wp(z),\wp'(z))\frac{dz}{\omega}$. Since $\wp(\omega-z)=\wp(-z)=\wp(z)$ and $\wp'(\omega-z)=\wp'(-z)=-\wp'(z)$, by a change of variables $z\to\omega-z$ we obtain $\lim_{N\to\infty}\sum_{n=1}^{N-1}f(\wp(z_n),\wp'(z_n))=\int_{\omega/2}^{\omega}f(\wp(z),\wp'(z))\frac{dz}{\omega}+f(\wp(z),-\wp'(z))\frac{dz}{\omega}$. Now since $x=\wp(z)$ and $\wp'(z)=y$ we have $dx=\wp'(z)dz \implies dz=\frac{dx}{y}$ and thus noting that $\wp(\omega/2)=e$ ($e=e_1$ or $e=e_2$ depending on the case of $\omega_{1,2}$) and $\wp(0)=O$ we get \\ $\lim_{N\to\infty}\sum_{n=1}^{N}f(\wp(z_n),\wp'(z_n))=\frac{1}{\omega}\int_{e}^{\infty}(f(x,y)+f(x,-y))\frac{dx}{y}$. The condition $\frac{f(x,y)}{y}\in{}\mathcal{IR}((e,\infty))$ and $f$ bounded is sufficient since $f(\wp,\wp')$ is bounded in $[0,\omega]$ iff $f$ is bounded in $[e,\infty)$ and for every for closed interval $[a,b]\subset(\frac{\omega}{2},\omega)$ we have $f(\wp,\wp')\in\mathcal{R}([a,b]) \iff \frac{f(x,y)}{y}\in\mathcal{R}([\wp(a),\wp(b)])$. This leaves only the problematic bounds $0,\omega$ where $y$ or $x$ is not bounded, where improper integration is still well defined however.
\end{proof}
\end{theorem}
\begin{observation}
Notice that $f$ can naturally be a complex valued function $f:[e,\infty)\to\mathbb{C}$ resulting in a complex integral over the real line.
\end{observation}
\begin{corollary}\label{conc}
In the particular case of $z_n=nz$ we have
\begin{equation} \label{eqimp}
    \lim_{N\to\infty}\frac{1}{N}\sum_{n=1}^{N}f(x_{nP},y_{nP})=\frac{1}{\omega}\int_{e}^{\infty}(f(x,y)+f(x,-y))\frac{dx}{y}, y\geq{}0
\end{equation}
\end{corollary}
Setting $f(x,y)=\one_{[0,\omega]}$ we get $1=\lim_{N\to\infty}\frac{1}{N}\sum_{n=1}^{N}1=\frac{1}{\omega}\int_{e}^{\infty}\frac{dx}{y} \implies \int_{e}^{\infty}\frac{dx}{y}=\omega$ (with $y$ taking values in
the whole $\mathbb{R}$) which is a result that is immediate by the Uniformization Theorem.

\begin{observation}\label{rem}
The sequence $z_n$ is equidistributed in $[0,\omega]$ iff $az_n$ is equidistributed in $[0,a\omega]$. For an elliptic curve $E_1$ with lattice $\Lambda_1$ every isomorphic elliptic curve is of the form $\Lambda_2=a\Lambda_1$. The isomorphism is the map $z\mod\Lambda_1\to{}az\mod{}\Lambda_2$ and so we get that: $P_n$ is equidistributed in $E_1$ w.r.t. the measure $\mu(X)=\frac{1}{\omega}\int_{\wp(X)}\frac{dx}{y_1}$ iff $\phi(P_n)$ is equidistributed in $E_2$ w.r.t. the measure $\mu(X)=\frac{1}{|a|\omega}\int_{\wp(X)}\frac{dx}{y_2}$.
\end{observation}

\subsection{Equidistribution in the whole space $E(\mathbb{R})$}

 We will now analyze the space $(A_1,\mathcal{B}_1,\mu_+)$ as defined in Theorem \ref{ththe}.

 \begin{theorem}
 The sequence $s_n=nz$ is equidistributed in $A_1$ iff $z\in\Gamma_1$ and $z\not\in\omega\mathbb{Q}$.
 \begin{proof}
 ($\implies$) This direction is obvious from Theorem \ref{ththe}.
 ($\impliedby$) Suppose and $z\in\Gamma_1$ and $z\not\in\omega\mathbb{Q}$. We observe that $s_{2n}\in{}B_1$ and $s_{2n+1}\in\Gamma_1$ and $s_{2n}=2s_n$, $s_{2n+1}=s_{2n}+z$. However $s_n$ equidistributed implies $ks_n$ is also equidistributed for every $k\in\mathbb{Z}$ and thus $s_{2n}$ is equidistributed in $B_1$ and $s_{2n+1}$ is equidistributed in $\Gamma_1$.
 \end{proof}
 \end{theorem}

 We then get the following theorem:
 \begin{theorem}
 Let $P\in{}E(\mathbb{R}): x_P\in{}(e_3,e_2)$ and $P\not\in{}E_{Tor}$, then\\ $\lim_{N\to\infty}\frac{1}{N}\sum_{n=1}^{N}f(x_{nP},y_{nP})=\frac{1}{2\omega}(\int_{e_3}^{e_2}(f(x,y)+f(x,-y))\frac{dx}{y}+\int_{e}^{\infty}(f(x,y)+f(x,-y))\frac{dx}{y}), \\ y\geq{}0$ for every bounded function $f\in\mathcal{R}((e_3,e_2)\cup{}(e,\infty))$.
 \end{theorem}

\subsection{Equidistribution in Curves}
The primary problem that arises here is that a curve may not be a probability space as it can be isomorphic to $\mathbb{R}$ in the topological sense with $\gamma(t)=(x_1(t),..,x_n(t))$ and $\lim_{t\to{}1}x_i(t)=\infty$ or $\lim_{t\to{}0}x_i(t)=\infty$. We may define curves on $\mathbb{P}^2$ in which case we have $\gamma(0)=O$ or $\gamma(1)=O$ as to be compliant with the definition of a curve but we will study them as affine curves through the natural map $(x_1(t),..,x_n(t),1)\to{}(x_1(t),..,x_n(t))$.
We will bypass the problem by defining equidistribution in a manner suitable for a non-compact space, like one isomorphic to $\mathbb{R}$ for example, in a manner similar to Gerl \cite{Gerl}.

\begin{definition} (\textbf{Gerl})\\
 Let $(X,\mathcal{B},\mu)$ be a measure space where $X$ is a locally compact Hausdorff space with countable base and $\mu$ a Radon measure (possibly not finite). Then a sequence $s_n$ is equidistributed in $X$ w.r.t. $\mu$ iff for every pair of compact subsets $A,B\subseteq{}X$ with $\mu(\partial{}A)=\mu(\partial{}B)=0$ we have
  \begin{equation}
 \lim_{n\to\infty}\frac{|\{s_1,\dots,s_n\}\cap{}B\cap{}A|}{|\{s_1,\dots,s_n\}\cap{}A|}=\mu(B)/\mu(A)
 \end{equation}
\end{definition}
 Since we are only interested in topological spaces like $\mathbb{R}$ and only need a definition for intervals of the form $[a,b]$ which always have trivial boundary $a,b$, we can use a more simple version. We can also drop the "for every pair of subsets" in favour of an increasing family of open intervals that covers the space since it will eventually contain any two such intervals. Before stating this definition we will define the problematic measure in the case of a curve:

\begin{theorem}
A continuous curve $\gamma:[0,1]\to\mathbb{R}^2$ equipped with the Borel $\sigma-$algebra of open sets of the curve is a measure space with respect to the Radon measure $\mu_{\gamma}(X)=\int_0^{1}||\gamma'(t)||\one_{X}(t)dt$. We thus obtain a measure space $(\gamma,\mathbb{B},\mu_{\gamma})$.
\begin{proof}
Obviously $\mu_{\gamma}(X)\geq{}0, \forall{}X\in\mathbb{B}$ and $\mu_{\gamma}(\varnothing)=0$. However for any countable collection of sets $\{A_i\}_1^{\infty}$ we obtain that  $\mu_{\gamma}(\cup_{k=1}^{\infty}A_k)=\int_0^{1}||\gamma'(t)||\sum_{k=1}^{\infty}\one_{A_k}(t)dt=\sum_{k=1}^{\infty}\mu_{\gamma}(A_k)$ where the interchange between the sum and the integral follows by Tonelli's Theorem since for the functions $f_n(t)=||\gamma'(t)||\one_{A_n}(t)$ we have $f_n\geq{}0$. The Radon property is obvious by the continuity of the curve.
\end{proof}
\end{theorem}

\begin{definition} \label{defin}
 A sequence of points $u_n=\gamma(s_n)$ defined on a curve $\gamma$ given by a sequence $s_n\in[0,1]$ is equidistributed iff
 \begin{equation} \label{equid}
 \lim_{n\to\infty}\frac{|\{s_1,\dots,s_n\}\cap{}[a,b]\cap{}A_k|}{|\{s_1,\dots,s_n\}\cap{}A_k|}=\mu_{\gamma}([a,b]\cap{}A_k)/\mu_{\gamma}(A_k)
 \end{equation}
 for every $[a,b]\subseteq[0,1]$ and every family of intervals $A_k=(a_k,b_k), a_k\ne{}0, b_k\ne{}1$, $A_k\subseteq{}A_{k+1}$ with $\cup_{k=1}^{\infty}A_k=(0,1)$.
\end{definition}

The information Definition \ref{defin} encodes is that every interval contains a proportion of the sequence proportionate to ``how much" of the curve is over that interval.

\begin{lemma}
 Definition \ref{defin} is not dependent on the set family $A_k$. More formally if $s_n$ is equidistributed w.r.t. $\mu_{\gamma}$ and a family of intervals $A_k$, then if $B_k$ is another family of intervals with the same properties, $s_n$ is also equidistributed w.r.t. $\mu_{\gamma}$ and $B_k$.
\begin{proof}
Suppose $A_k,B_k$ are two such families, then $\mu(A_k)<\mu([0,1])=1$ and $\lim_{k\to\infty}\mu(A_k)=\mu(\cup_{k=1}^{\infty}A_k)=1$ which implies that $\forall{}\epsilon\in{}(0,1), \exists{}k\in\mathbb{N}: \mu(A_k)\in{}(1-\epsilon,1)$. So for every $B_i=(a_i,b_i)$ setting $\epsilon=min\{a_i-0,1-b_i\}$ there exists a $k\in\mathbb{N}: \mu(A_k)\geq{}1-\epsilon/2$ and thus if $A_k=(c_k,d_k))$ and $c_k\geq{}a_i$ or $d_k\leq{}b_i$ we would have $1-\epsilon/2\leq{}\mu(A_k)\leq{}1-\epsilon \implies e/2\leq{}0$ contradiction. Thus $B_i\subset{}A_k$ and so supposing Equation \ref{equid} holds for every $A_k$ it also hold for all $B_i$. The same argument for $B_k$ completes the proof.
\end{proof}
\end{lemma}

We would like to emphasize how this definition is a natural extension of the definition of equidistribution for a compact space since in that case we obtain the usual definition by setting $A_k$ equal to our space.
With the above lemma we can choose symmetric $A_k$ that will make integration easier on the real line. We will thus only consider families of intervals $A_k=(\omega/2-a_k,\omega/2+a_k)$ with $a_k$ increasing and $a_k<\omega/2$ thus attaining $\lim_{k\to\infty}a_k=\omega/2$.

We now take a look at an example which showcases what happens a sequence equidistributed in $\mathbb{R}/\mathbb{Z}$ when projected on a circle. The following example is what motivated the use of the $\mu_{\gamma}$ measure in our definition:
\begin{example}
Let $\gamma(t)=(sin(2\pi{}t)$ and $cos(2\pi{}t))$ $g:[0,1]\to{}S^1$. Then for any sequence $t_n$ equidistributed in $\mathbb{R}/\mathbb{Z}$ and $f\in\mathcal{R}([0,1])$ we have $\lim_{N}\sum_{n=1}^Nf(t_n)=\int_0^1f(\gamma(t))dt$. Notice that setting $sin(2\pi{}t)=x, cos(2\pi{}t)=y \implies dx/y=dx/cos(2\pi{}t)=dt$ and thus the integral becomes (integrating along $y>0$ and $y<0$ as before) $\int_{-1}^1f(x,y)+f(x,-y)\frac{dx}{y}$. We make the following observation: $||\gamma'(t)||=\sqrt{1+(\frac{dy}{dx})^2}=\frac{1}{y}\sqrt{x^2+y^2}=\frac{1}{y}, y>0$. Indeed then we obtain the expected formula \\ $\lim_{N\to\infty}\sum_{n=1}^Nf(t_n)=\int_0^1fd\mu_{\gamma}$.
\end{example}

\begin{theorem} \label{lol}
A sequence of points $u_n=\gamma(s_n)$ defined on a curve $\gamma$ given by a sequence $s_n\in[0,1]$ is equidistributed iff $\lim_{N\to\infty}\frac{\sum_{n=0}^{N-1}f_k(s_n)}{|\{s_1,\dots,s_n\}\cap{}A_k|}=\frac{\int_{A_k}fd\mu_{\gamma}}{\mu_{\gamma}(A_k)}$ where $f_k=f\one_{A_k}$ with $f$ Riemann integrable in every $A_k$.
\begin{proof}
We observe that setting $E=\mu_{\gamma}(A_k)$ we obtain a probability space \\ $(A_k,\mathcal{B}\cap{}A_k,\mu/E)$ and the result is then an immediate consequence of Weyl's Theorem.
\end{proof}
\end{theorem}

\subsection{Equidistribution in Real Elliptic Curves}

In the case of an Elliptic Curve, Theorem \ref{lol} is phrased as:
\begin{corollary}
 A sequence of points $u_n=(\wp(s_n),\wp'(s_n))$ defined on an Elliptic Curve given by a sequence $s_n\in[0,\omega]$ is equidistributed iff
 \begin{equation} \label{eq1}
\lim_{N\to\infty}\frac{\sum_{n=0}^{N-1}f_k(u_n)}{|\{s_1,\dots,s_n\}\cap{}A_k|}=\frac{\int_e^{\wp(a_k)}(f(x,y)+f(x,-y))\sqrt{1+(\frac{dy}{dx})^2}dx}{2\int_e^{\wp(a_k)}\sqrt{1+(\frac{dy}{dx})^2}dx}
\end{equation} where $y\geq{}0$ and $f\in\mathcal{IR}((e,\wp(a_k)))$ and bounded, $\forall{}k\in\mathbb{N}$.
\end{corollary}

We then have from Equation (\ref{eqimp}) that
\begin{equation} \label{eq2}
\lim_{N\to\infty}\frac{\sum_{n=0}^{N-1}f_k(u_n)}{|\{s_1,\dots,s_n\}\cap{}A_k|}=\frac{\int_e^{\wp(a_k)}f(x,y)+f(x,-y)\frac{dx}{y}}{2\int_e^{\wp(a_k)}\frac{dx}{y}}, y\geq{}0
\end{equation}

\begin{theorem}
The points of the sequence $s_n=[n]P$ where $P\in{}E$, $E$ an Elliptic Curve are not equidistributed on $E$ with respect to the "natural" measure $\mu_{\gamma}$ but are instead equidistributed with respect to the measure $\mu(X)=\frac{1}{\omega}\int_{\wp(X)}\frac{dx}{y}$.
\begin{proof}
The result follows from an immediate comparison of Equations \ref{eq1}, \ref{eq2}. Choosing $f$ as the indicator function of some interval $[a,b]$ and taking the limit $\lim_{k\to\infty}a_k$ in both cases we get $0$ from Equation \ref{eq1} and $\frac{\mu(\wp^{-1}([a,b]))}{\omega}>0$ from Equation \ref{eq2}.
\end{proof}
\end{theorem}

Indeed the points of $[n]P$ are tightly concentrated around $e$ and get thinner and thinner as we approach infinity. However the sequence remains dense in every set $[a,b]\subseteq{}[e,\infty)$.

A new question arises now: Can we possibly equip $[0,\omega]$ with a different measure $\mu'$ such that
$\lim_{N\to\infty}\frac{1}{N}\sum_{n=1}^{\infty}f_k(\wp(nz))=\frac{\int_{A_k}f\circ\wp{}d\mu'}{\mu'(A_k)}$ for every $k\in\mathbb{N}$?\\
In the case of probability measures the answer is negative since a sequence in a compact space is equidistributed w.r.t. at most one probability measure. To see this we note that $\mu((a,b))=\mu([a,b])=\mu'((a,b))$ for every open set in $[0,\omega]$ and open sets generate the Borel $\sigma-$algebra which is stable under finite intersection. As a consequence of the monotone class theorem the measures $\mu,\mu'$ agree on every set of $\mathcal{B}$.
Observe that by the Riesz Representation Theorem, changing the measure is equivalent to sampling by a different positive function since $\int_0^{\omega}f\circ{}gd\mu=\int_0^{\omega}d\mu_{g}$ as a positive linear functional.
Even in the case of a Radon measure we get the following:
\begin{theorem} \label{th4}
Let $P_n$ be equidistributed in $E(\mathbb{R})$ w.r.t. $\mu(X)=\frac{1}{\omega}\int_{\wp(X)}\frac{dx}{y}$, then there exists no function $f$ with $\mu(\partial{}f^{-1}(A))=0$ for all compact intervals $A\subset{}[e,\infty)$ taking $x_n\to{}f(x_n)$, such that $f(x_n)$ is equidistributed w.r.t. any non finite, Radon measure $\mu_R$.
\begin{proof}
Suppose such a function exists. Then $\one_{f^{-1}(A)}/y$ is improper Riemann integrable in $(e,\infty)$
We first observe that for any closed intervals $A\subseteq{}B\subset{}[e,\infty)$ it holds that: \\ $\lim_{N\to\infty}\frac{|\{f(x_1),\dots,f(x_N)\}\cap{}A|}{|\{f(x_1),\dots,f(x_N)\}\cap{}B|}=\lim_{N\to\infty}\frac{\frac{1}{N}\sum_{n=1}^N\one_{A}(f(x_n))}{\frac{1}{N}\sum_{n=1}^N\one_{B}(f(x_n))}=\frac{\int_{f^{-1}(A)}\frac{dx}{y}}{\int_{f^{-1}(B)}\frac{dx}{y}}$. Now by Definition \ref{defin} it must be the case that $\frac{\int_{f^{-1}(A)}\frac{dx}{y}}{\int_{f^{-1}(B)}\frac{dx}{y}}=\frac{\mu_R(A)}{\mu_R(B)}$. Taking a sequence $B_n\to[e,\infty)$ we must then have that $\lim_{n\to\infty}\frac{\int_{f^{-1}(A)}\frac{dx}{y}}{\int_{f^{-1}(B_n)}\frac{dx}{y}}=0$. However $f^{-1}(B_n)\subset{}[e,\infty)$ and thus $\int_{f^{-1}(B_n)}\frac{dx}{y}<\omega$ contradicting our previous claim.
\end{proof}
\end{theorem}
A function that would contradict Theorem \ref{th4} would obviously satisfy $\mu(\partial{}f^{-1}(A))>0$ for some closed interval $A$ and thus since $\partial{}f^{-1}(A)\subseteq{}f^{-1}(\partial{}A)$ we must have $\mu(f^{-1}(\partial{}A))>0$. This implies that if $A=[a_1,a_2]$, then the set $U=\{x>e:f(x)=a_1\text{ or }f(x)=a_2\}$ has positive measure. Thus $f$ is clearly either discontinuous in a positive measure subset of points or $f$ changes monotonicity in a positive measure subset of points or $f$ is nowhere monotonic. This aims to show that $f$ is not trivial to find.

\subsection{Distribution of Polynomial Maps on Elliptic Curves}

All of our previous theorems are phrased for an equidistributed sequence in $[0,\omega]$ in general. This enables our previous theorems to be restated for any polynomial sequence on an elliptic curve:
\begin{theorem}(\textbf{Weyl's Equidistribution Theorem})\\
Let $p(x)$ be a monic polynomial in $\mathbb{Z}[x]$, then the sequence $p(n)u$ is equidistributed in $[0,1]$ iff $u\not\in\mathbb{Q}$.
\end{theorem}

\begin{theorem}
Let $p(x)$ be a monic polynomial in $\mathbb{Z}[x]$, then the sequence $p(n)z$ is equidistributed in $[0,\omega]$ iff $z\not\in\omega\mathbb{Q}$.
\begin{proof}
The proof is an immediate modification of the original proof in the case of $[0,1]$. For the full proof see Corollary 3 of \cite{ref_url}.
\end{proof}
\end{theorem}
This means that every for every monic polynomial in $\mathbb{Z}$ we have the following corollary:
\begin{corollary}
 Let $p(x)$ be a monic polynomial in $\mathbb{Z}[x]$ then the sequence $s_n=[p(n)]P$ is equidistributed w.r.t. $\mu(X)=\frac{1}{\omega}\int_{\wp(X)}\frac{dx}{y}$ iff $\wp(P)\not\in\omega\mathbb{Q}$.
\end{corollary}

\subsection{Equidistribution of points in $E(\mathbb{Q})$}
When working with computers there is an obvious limitation to the field of rationals $\mathbb{Q}$. This actually makes things easier since we can specifically state which points will give equidistributed sequences in $E(\mathbb{Q})$ with respect to the measure $\mu(X)=\frac{1}{\omega}\int_{X}\frac{dx}{y}$. Let us make clear something ambiguous first:

\begin{definition}
 We say that a sequence $s_n\in{}E(\mathbb{Q})$ is equidistributed in $E(\mathbb{Q})$ w.r.t. a measure $\mu$ iff $s_n$ is equidistributed in $E(\mathbb{R})$ w.r.t. the measure $\mu$.
\end{definition}
Thus restricted to $Q$ we use the dynamics of it's extension $\mathbb{R}$ to define equidistribution for our purposes.
By the Mordell-Weil Theorem (page 220 of \cite{silverman}) we know that $E(\mathbb{Q})=E_{Tor}\oplus\mathbb{Z}^r$, and so:
\begin{theorem}
A point $P\in{}E(\mathbb{Q})$ is equidistributed w.r.t $\mu$ in $E(\mathbb{R})$ iff $P\not\in{}E(\mathbb{Q})_{Tor}$. Thus $\forall{}P\in{}E(\mathbb{Q})$ with $y_P\ne{}0$:
\begin{itemize}
    \item $x_P,y_P\not\in\mathbb{Z}$ or
    \item $x_P,y_P\in\mathbb{Z}$ but $y^2\nmid\Delta$
\end{itemize}
the sequence $\psi_n=nP$ is equidistributed w.r.t $\mu$.
\begin{proof}
By Theorem \ref{th2} we have that if $P=\wp(z_0)$ then $\psi_n$ is equidistributed in $E(\mathbb{R})$ iff $z_0\not\in\omega\mathbb{Q}$. We now observe that $z_0\in\omega\mathbb{Q} \iff nz_0=0\mod{}[0,\omega] \iff nP=O \iff P\in{}E(\mathbb{Q})_{Tor}$. An immediate application of Nagell-Lutz now completes the theorem.
\end{proof}
\end{theorem}

\section{Distribution in $\mathbb{R}/\mathbb{Z}$} \label{sec}

Suppose $x_n=\wp^{-1}(s_n)$ where $s_n$ is equidistributed in $[0,\omega]$. We will investigate if such a sequence could produce a sufficiently good PRNG $\mod(1)$. Let us first examine the most simple case of a sequence taking $E(\mathbb{R})\to\mathbb{R}/\mathbb{Z}$: $\phi_n=\{x_n\}$.
By Weyl's Criterion for equidistribution we want to show that: $\lim_{N\to\infty}\frac{1}{N}\sum_{n=1}^{N}e^{2\pi{}jkx_{n}}=0, \forall{}k\in\mathbb{Z}$. By Equation \ref{eqimp} we then need to show that $\int_e^{\infty}\frac{e^{2\pi{}jkx}}{y}dx=0 \iff \int_0^{\infty}\frac{e^{2\pi{}jkx}}{\sqrt{x(x+a)(x+b))}}dx=0$ by a simple change of variables $x\to{}x+e$. We see that $a,b<0$ since $0$ is now the largest root of $y^2=0$. We observe however that this cannot happen when $y>0$ is increasing since integrating over a period $\Im(\int_n^{n+1}\frac{e^{2\pi{}jkx}}{\sqrt{x(x+a)(x+b)}}dx)>0, \forall{}n\in\mathbb{Z}$ implying $\Im(\int_0^{\infty}\frac{e^{2\pi{}jkx}}{\sqrt{x(x+a)(x+b)}}dx)=\\ \Im(\sum_{n=1}^{\infty}\int_n^{n+1}\frac{e^{2\pi{}jkx}}{\sqrt{x(x+a)(x+b)}}dx)>0$.
\\ With this in mind we seperate two cases:
\begin{lemma} \label{ellitpic}
Let $y(x)=\sqrt{(x-e)(x-e_1)(x-e_2)}, e>e_1>e_2$ or $e>0$ and $e_1,e_2\not\in\mathbb{R}$ and $y^2=x^3+Ax+B$. Then:
\begin{itemize}
    \item $y$ is increasing iff $A\geq{}0$ or $\sqrt{\frac{|A|}{3}}<e$ when $A<0$.
    \item $y$ is increasing in $[e,-\sqrt{\frac{|A|}{3}})\cup{}[\sqrt{\frac{|A|}{3}},\infty)$ and decreasing in \\ $(-\sqrt{\frac{|A|}{3}},\sqrt{\frac{|A|}{3}})$ otherwise
\end{itemize}
\begin{proof}
It is obvious that $\frac{dy}{dx}=\frac{3x^2+A}{2y}$ and thus $y$ is decreasing in $[-\sqrt{\frac{|A|}{3}},\sqrt{\frac{|A|}{3}})$. Thus for $y$ is increasing iff $[e,\infty)\cap{}[-\sqrt{\frac{|A|}{3}},\sqrt{\frac{|A|}{3}})=\emptyset$. Next notice that $e\not\in[-\sqrt{\frac{|A|}{3}},\sqrt{\frac{|A|}{3}})$, otherwise $y$ would be decreasing in $[e,\sqrt{\frac{|A|}{3}})$ and then $\forall{}x\in{}[e,\sqrt{\frac{|A|}{3}})$ we would have $y(x)<0$. So the above are indeed the only two cases.
\end{proof}
\end{lemma}

 \begin{figure}[htp]
\centering
\includegraphics[width=4cm]{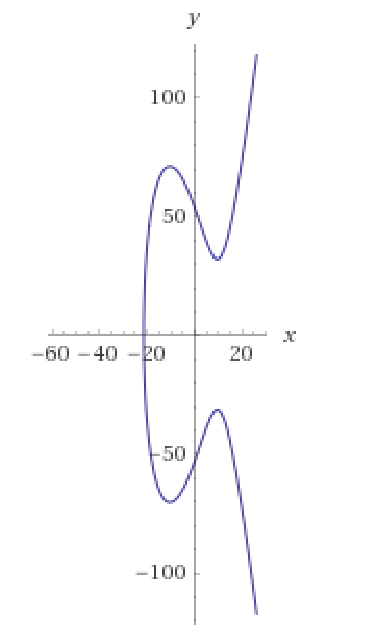}
\caption{An elliptic curve where $[e,\infty)\cap{}[-\sqrt{\frac{|A|}{3}},\sqrt{\frac{|A|}{3}})\ne{}\emptyset$}
\label{fig:ellitpic}
\end{figure}

\begin{lemma} \label{pol}
Suppose $s_n=\wp^{-1}(P_n)$ (where $P_n=(x_n,y_n,1)$) is equidistributed in $[0,\omega]$ as defined by an elliptic curve $y^2=x^3+Ax+B$ with $A\geq{}0$ or $\sqrt{\frac{|A|}{3}}<e$, then the sequence $\{x_n\}$ is not equidistributed in $\mathbb{R}/\mathbb{Z}$ w.r.t the Lebesgue measure.
\begin{proof}
We basicaly restate what was written above in more general context. By Lemma \ref{ellitpic} we have that $y$ is increasing in $[e,\infty)$.
\\By Weyl's criterion $x_n$ is equidistributed iff $\lim_{N\to\infty}\frac{1}{N}\sum_{n=1}^{N}e^{2\pi{}jkx_{n}}=0, \forall{}k\in\mathbb{Z}$. From Equation \ref{eqimp} however it follows that \\ $\lim_{N\to\infty}\frac{1}{N}|\sum_{n=1}^{N}e^{2\pi{}jkx_{n}}|=\frac{2}{\omega}|\int_e^{\infty}\frac{e^{2\pi{}jkx}}{y}dx|>0$. It remains to prove that $\Im(\int_n^{n+1}\frac{e^{2\pi{}jkx}}{\sqrt{x(x+a)(x+b)}}dx)>0$. This is equivalent to showing that \\ $\int_n^{n+1}\frac{sin(2\pi{}kx)}{\sqrt{x(x+a)(x+b)}}dx>0$. We observe that: 
$\int_n^{n+1}\frac{sin(2\pi{}kx)}{\sqrt{x(x+a)(x+b))}}dx=\\ \int_n^{n+\frac{1}{2}}\frac{sin(2\pi{}kx)}{\sqrt{x(x+a)(x+b)}}dx+\int_{n+\frac{1}{2}}^{n+1}\frac{sin(2\pi{}kx)}{\sqrt{x(x+a)(x+b)}}dx= \int_n^{n+\frac{1}{2}}sin(2\pi{}kx)(\frac{1}{\sqrt{x(x+a)(x+b)}}-\frac{1}{\sqrt{(x+\frac{1}{2})(x+a+\frac{1}{2})(x+b+\frac{1}{2})}}dx)>0$ since $y(x)<y(x+\frac{1}{2})$.
\end{proof}
\end{lemma}

We can now pass to the case of three distinct real roots:
\begin{theorem} \label{thth}
Let $s_n=\wp^{-1}(P_n)$ (where $P_n=(x_n,y_n,1)$) be equidistributed in $[0,\omega]$ defined by an elliptic curve $E: y^2=x^3+Ax+B$. Then if $E$ has $3$ distinct real roots, the sequence $\{x_n\}$ is not equidistributed w.r.t. the Lebesque measure.
\begin{proof}
Considering the function $f(x)=x^3+Ax+B$ the only way for it to have three real roots $e_1>e_2>e_3$ is iff $e_1>\sqrt{\frac{|A|}{3}}, e_2\in{}(-\sqrt{\frac{|A|}{3}},\sqrt{\frac{|A|}{3}})$ and $e_3<-\sqrt{\frac{|A|}{3}}$. The result is now obvious from Lemma \ref{pol}.
\end{proof}
\end{theorem}

An immediate indication of this result is the following:
\begin{observation}
Since $\lim_{N\to\infty}\frac{1}{N}\sum_{n=1}^N\one_{[a,b]}(\{x_n\})=\frac{1}{\omega}\sum_{n=\lfloor e \rfloor}^{\infty}\int_{n+a}^{n+b}
\frac{dx}{y}$. Defining the function $F_n(t)=\int_{n+a+t}^{n+b+t}\frac{dx}{y}$ where $t\in [b,1]$ we have that ${F_n}'(t)=\frac{1}{y(a+n+t)}-\frac{1}{y(b+n+t)}<0$ since $y$ is increasing. So $F_n$ is decreasing. Then choosing $a=0,b=\frac{1}{2}$ and $t=\frac{1}{2}$ gives $\int_{n}^{n+\frac{1}{2}}\frac{dx}{y}>\int_{n+\frac{1}{2}}^{n+1}\frac{dx}{y}$ which implies \\
$\lim_{N\to\infty}\frac{1}{N}\sum_{n=1}^N\one_{[0,1/2]}(\{x_n\})>\lim_{N\to\infty}\frac{1}{N}\sum_{n=1}^N\one_{[1/2,1]}(\{x_n\})$.
\end{observation}

This still leaves us to deal with the case $\Delta<0$. This situation is much more complicated since we can't use the monotonicity of $y$. We will attenmpt a different approach.
\begin{lemma} \label{lem2}
Suppose $\int_e^{\infty}\frac{e^{2\pi{}jnx}}{y}dx=0$ for every $n\in\mathbb{Z}$. Then $\int_e^{\infty}\frac{f(x)}{y}dx=\frac{\omega}{2}\int_0^1f(x)dx$, for every $f\in\mathcal{L}^1([0,1])$.
\begin{proof}
Since trigonometric polynomials are dense in $\mathcal{L}^1([0,1])$ we have that for every $\epsilon>0$, then there exists a trigonometric polynomial \\ $p_N(x)=\sum_{n=-N}^Nd_{N,n}x^n$ such that $|p_N(e^{2\pi{}jx})-f(x)|<\epsilon$. By integrating we obtain $|\int_0^1f(x)dx-d_{N,0}|<\epsilon$. Dividing by $y>0$ and integrating we get $|\int_e^{\infty}\frac{f(x)}{y}dx-\frac{\omega}{2}d_{N,0}|<\frac{\omega}{2}\epsilon$ and finally with the triangle inequality: $|\int_e^{\infty}\frac{f(x)}{y}dx-\frac{\omega}{2}\int_0^1f(x)dx|<|\int_e^{\infty}\frac{f(x)}{y}dx-\frac{\omega}{2}d_{N,0}|+|\frac{\omega}{2}\int_0^1f(x)dx-\frac{\omega}{2}d_{N,0}|<\omega\epsilon$ and since $\epsilon$ is arbitary, the proof is complete.
\end{proof}
\end{lemma}

\begin{theorem} {(\textbf{General Version})} \label{genver}
Let $s_n=\wp^{-1}(P_n)$ (where $P_n=(x_n,y_n,1)$) be equidistributed in $[0,\omega]$. Then the sequence $\{x_n\}$ is not equidistributed in $[0,1]$ w.r.t. the Lebesque measure.
\begin{proof}
The case where $y$ is increasing is settled by Theorem \ref{thth}.
Otherwise suppose $P_n$ is equidistributed in $E$ and $\{x_n\}$ is also equidistributed in $[0,1]$. Define the map $E\to{}E':\phi(x,y)=(u^2x,u^3y)$, then $\phi$ is an isomorphism between $E$ and $E':y^2=x^3+\frac{A}{u^4}x+\frac{B}{u^6}$ for every $u>0$. Thus by Remark \ref{rem}, the sequence $\phi(P_n)$ is also equidistributed in $E'$ w.r.t. the measure $\mu(X)=\frac{1}{\omega}\int_{\wp(X)}\frac{dx}{y}$. Notice that for $u^2\in\mathbb{Z}$, if $\{x_n\}$ is equidistributed in $[0,1]$ we have that $\{u^2x_n\}$ is also equidistributed in $[0,1]$. We thus focus our attention to equidistributed sequences on $E'$. Then $y$ is increasing in $[e',-\frac{1}{u^2}\sqrt{\frac{A}{3}}]$,$[\frac{1}{u^2}\sqrt{\frac{A}{3}},+\infty)$ and decreasing in $(-\frac{1}{u^2}\sqrt{\frac{|A|}{3}},\frac{1}{u^2}\sqrt{\frac{|A|}{3}})$  where $e'$ is the largest root of $E'$. After centering the curve as before by setting $x\to{}x+e'$ (since $\int_{e'}^{\infty}\frac{e^{2\pi{}jnx}}{y}dx=e^{2\pi{}jne'}\int_0^{\infty}\frac{e^{2\pi{}jnx}}{y(x+e')}dx$), we have that $y$ is decreasing in $(-e'-\frac{1}{u^2}\sqrt{\frac{|A|}{3}},-e'+\frac{1}{u^2}\sqrt{\frac{|A|}{3}})$. Now pick an integer $u$ large enough so that $2\frac{1}{u^2}\sqrt{\frac{|A|}{3}}<\frac{1}{3}$. Partitioning $[0,1]$ into three equal length intervals we obtain that there is at least one interval $U$ such that $(-e'-\frac{1}{u^2}\sqrt{\frac{|A|}{3}},-e'+\frac{1}{u^2}\sqrt{\frac{|A|}{3}})\cap{}U=\emptyset$. However by Lemma \ref{lem2} for any periodic function $f:[0,1]\to\mathbb{R}$ with $\int_0^1f(x)dx=0$ we have $\int_{e'}^{\infty}\frac{f(x)}{y}dx=0$. We define the function $f(x)=\one_{U_1}(x)-\one_{U_2}(x)$ where $U_1=[a,\frac{a+c}{2})$ (the "right half"), $U_2=(\frac{a+c}{2},c]$ (the "left half") and $U=[a,c]$. We then expand $f$ to a function on $\mathbb{R}$ as $f_{\mathbb{R}}(x)=\one_{U_1}(\{x\})-\one_{U_2}(\{x\})$. Now as before (centering the curve at $0$ for convenience)  $\int_0^{\infty}\frac{f_{\mathbb{R}}(x)}{y}dx=\sum_{n=0}^{\infty}\int_{n}^{n+1}\frac{f(x)dx}{y}=\sum_{n=0}^{\infty}(\int_{U_1}\frac{dx}{y}-\int_{U_2}\frac{dx}{y})$ and since $(-e'-\frac{1}{u^2}\sqrt{\frac{|A|}{3}},-e'+\frac{1}{u^2}\sqrt{\frac{|A|}{3}})\cap{}(U_1\cup{}U_2)=\emptyset$ we always have that $\int_{n+U_1}\frac{dx}{y}-\int_{n+U_2}\frac{dx}{y})>0$. This implies that $\{u^2x_n\}$ is not equidistributed in $[0,1]$ and thus neither is $\{x_n\}$.
\end{proof}
\end{theorem}

Another possible question now is the following: Can we "fix" this sequence by taking the least significant digits that should exhibit more "random" behaviour?
The answer to that question is "no" since in that case we would essentialy require \\ $\lim_{N\to\infty}\frac{1}{N}\sum_{n=1}^{N}e^{2\pi{}jkx_{n}10^m}=0$, $\forall{}k\in\mathbb{Z}$ which would then be equivalent to showing that $10^{3m}\int_{10^me}^{\infty}\frac{e^{2\pi{}jkx}}{x^3+10^{2m}Ax+10^{3m}B}dx=0$ which is the same as proving that an equidistributed sequence $u_n$ with respect to the measure \\ $\mu(X)=\frac{1}{\omega'}\int_X\frac{dx}{x^3+10^{2m}Ax+10^{3m}B}$ is equidistributed in $\mathbb{R}/\mathbb{Z}$.

\section{Conclusion}
After providing the conditions for equidistribution of $[m]P$ over $E(\mathbb{C})$ in terms of linear independance over $\mathbb{Z}$ we turned to the much more interesting case of $E(\mathbb{R})$. Here we obtained the main result stated in Corollary \ref{conc} and concluded that the points of $[m]P$ and any other sequence that is equidistributed on the borders of the complex lattice follow the distribution described by Equation \ref{eq2}. The generalization to polynomial sequences is immediate from Weyl's well known result. Finally Theorem \ref{genver} provides a further result on the distribution of the rational part of $x_{[n]P}$, namely that it is not equidistributed with respect to the Lebesgue measure on $[0,1]$.

\end{document}